\documentclass{article}
\usepackage[centertags]{amsmath}
\usepackage{amsfonts}
\usepackage{amssymb}
\usepackage{amsthm}
\usepackage{newlfont}

\usepackage{graphicx}
\usepackage{graphicx}
\usepackage{amsfonts}
\usepackage{amssymb}
\usepackage{amsthm}
\usepackage{xcolor}
\usepackage{amsmath}

\usepackage[english]{babel}

\newenvironment{pf}{\noindent{\bf Proof.}\enspace}{
\hfill$\Box$\medskip}

 \newtheorem{theorem}{Theorem}[section]
 \newtheorem{cor}[theorem]{Corollary}
 \newtheorem{lemma}[theorem]{Lemma}

 \theoremstyle{definition}
 \newtheorem{defn}{Definition}[section]
\newtheorem{remark}[theorem]{Remark}
\newtheorem{example}[theorem]{Example}

 \numberwithin{equation}{section}

  \newcommand{\cqfd}
{%

\mbox{}%
 \nolinebreak%
 \hfill%
 \rule{2mm}{2mm}%
 \par%
 }

\begin{document}

\title{\bf  Dynamics and oscillations of
models for \\
differential equations with delays}

\author{{\bf Mohsen Miraoui}$^{a,b}$ {\bf and  Du\v{s}an D. Repov\v{s}}$^{c,d,e}$\\
 $^a$ {\small IPEIK, University of Kairouan, Tunisia}\\
$^b${\small LR11ES53, FSS, University of Sfax, Tunisia}\\
{\small {\it miraoui.mohsen@yahoo.fr}}\\
$^c$ {\small Faculty of Education, University of Ljubljana,  Slovenia}\\
$^d$ {\small Faculty of Mathematics and Physics, University of Ljubljana, Slovenia}\\
$^e$ {\small Institute of Mathematics, Physics and Mechanics, Ljubljana, Slovenia}\\
{\small
 {\it dusan.repovs@guest.arnes.si}}
}

 \date{}

\maketitle

\begin{abstract}
By developing new efficient techniques and using an appropriate fixed point theorem, we derive several new sufficient conditions for the  pseudo almost periodic solutions with double measure for some system of differential equations with delays. As an application, we consider certain models for neural networks with delays.
\end{abstract}

\medskip


 
\section{Introduction}
Existence of periodic, almost periodic and pseudo almost periodic solutions of differential equations has great significance and is therefore an important problem. Such dynamics can be found in electronic
circuits  and many other physical and biological systems (see \cite{BM1,Boch,C,GMN, Ni,D2,D1,Ra,Z}). Ezzinbi et al. \cite{EM} introduced a new and powerful measure-theoretic method to resolve this open problem. Since then, this method has been used for various classes of evolution equations as well as stochastic differential equations
and has become very popular.

The notion of measure pseudo almost periodicity was first introduced by Blot et al. \cite{EM}
(see also \cite{Ait,F1,E2015,E2016,mir3,mir2,mir4, Zhao}). Obviously,
these new results generalize the earlier work od Diagana \cite{BM}. Recently, Diagana et al. \cite{mir} have  introduced the notion of
double measure pseudo almost periodicity  as a generalization of the
measure pseudo almost periodicity. We note that this generalized concept coincides with the latter one (take $\mu\equiv\nu$). 

 In this paper, by applying an appropriate fixed point theorem, we derive
  some
conditions which
ensure  the existence,  the exponential
stability, and  the uniqueness  of $(\mu,\nu)-$pap solutions of the following models  with delays:

\begin{eqnarray}\label{sfw}
  &&x_{i}'(t)= -c_{i}(t) x_{i}(t) + \sum_{j=1}^{n} d_{ij} (t) f_{j} (t,x_{j}(t)) +
\sum_{j=1}^{n} a_{ij}(t) g_{j} (t,x_{j} (t- \tau_{ij}) ) \nonumber\\
  &&+  \sum_{j=1}^{n} \sum_{l=1}^{n} b_{ijl}(t) h_{j}(t,x_{j}(t-
\sigma_{ij})) h_{l} (t,x_{l}(t- \nu_{ij})) + I_{i}(t),     \\
&&x_i(s)= \varphi_{i} (s) , \; s \in (- \theta , 0 ] , \   i\in \{1,...,n\},  \nonumber\end{eqnarray}
 where  functions
 $$c_i, I_i, d_{ij}, a_{ij},b_{ijl}:\mathbb{R}\to \mathbb{R}\ \mbox{ and} \
   f_j,g_j,h_j:\mathbb{R}\times \mathbb{R}\to \mathbb{R}, \ \ i,j,l\in \{1,...,n\}    $$
   are continuous and
    $\tau_{ij}$, $\sigma_{ij}$, and $\nu_{ij}$ are positive constants.

The paper is organized as follows: in Section 2 we collect key definitions, examples, and basic results. In Section 3 we discuss the existence, the stability and the uniqueness  of
double measure pseudo almost periodic
solutions of system \eqref{sfw}. Finally, in Section 4
we  present an application which  illustrates the effectiveness of our
results.

\section{ Preliminaries}

\begin{defn}(see \cite{EM})
Let $f$ be a continuous function on $\mathbb{R}$ with  values in $\mathbb{R}^n $.
Then  $f$ is said to be
{\it almost periodic}, denoted by $f \in \mathcal{AP}(\mathbb{R},\mathbb{R}^n)$, if for all $\varepsilon>0$,  there exists
a
 number
$ l(\varepsilon)>0$ such that
every  interval $I$ of length $l(\varepsilon)$
 contains a point  $\tau\in \mathbb{R}
$ with the property that  $$\| f(t + \tau) - f(t) \| < \varepsilon, \ \ \mbox{ for all} \ \
t
\in \mathbb{R}.$$
\end{defn}
The space  $\mathcal{AP}(\mathbb{R},\mathbb{R}^n)$ equipped with the norm
\begin{equation*}
\|f\|_{\infty}:=\max_{1\leq i \leq n}\sup_{t\in \mathbb{R}}|f_i(t)|
\end{equation*}
is then a Banach space.
Let $\mathcal{B}$ be
the Lebesque $\sigma$-field on
 $\mathbb{R}$ and define a collection $\mathcal{M}$ of measures on $\mathcal{B}$
$$\mathcal{M}=\{\mu \ \ \mbox{is a positive measure on}\;
\mathcal{B}; \ \mu(\mathbb{R})=+\infty, \ \
\mbox{and} \ \
\mu([s,t])<\infty,\;
 \mbox{for all} \ \;s,t \in \mathbb{R},\;s\leq t\}.$$

Let $X$ be a Banach space and denote by
$\mathcal{BC}(\mathbb{R},X)$
 the Banach space of bounded continuous functions from
 $\mathbb{R}$ to $X$, equipped with the supremum norm
$\| f \|_{\infty} = \sup_{t \in \mathbb{R}} \| f(t) \|$.
In order to be able to introduce   double measure pseudo almost periodic functions, we need  the following ergodic  spaces

$$\mathcal{E} ( \mathbb{R}, \mathbb{R}^{n}, \mu ,\nu) := \{ f \in \mathcal{BC} (\mathbb{R}, \mathbb{R}^{n}): \lim_{z \rightarrow \infty} \frac{1}{\nu([-z,z])} \int_{-z}^{z} \| f (t) \| d \mu(t) = 0 \}, $$
and
$$\mathcal{E} ( \mathbb{R}, \mathbb{R}^{n}, \mu ) := \mathcal{E} ( \mathbb{R}, \mathbb{R}^{n}, \mu ,\mu)=\{ f \in \mathcal{BC} (\mathbb{R}, \mathbb{R}^{n}): \lim_{z \rightarrow \infty} \frac{1}{\mu([-z,z])} \int_{-z}^{z} \| f (t) \| d \mu(t) = 0 \}. $$

\begin{defn}(see \cite{mir})
If $ \mu ,\nu\;\in \mathcal{M}$, then  $f \in
\mathcal{BC}(\mathbb{R},\mathbb{R}^{n})$  is said to be  $(\mu,\nu)-${\it pseudo almost periodic}, abbreviated as
$(\mu,\nu)-${\it pap},
denoted by $f\in \mathcal{PAP} (\mathbb{R},\mathbb{R}^{n}, \mu,\nu)$,
if there exists a decomposition
\begin{equation}\label{unique}
  f = g + \varphi, \
  \mbox{where} \
   \varphi \in
\mathcal{E} ( \mathbb{R},\mathbb{R}^{n}, \mu,\nu) \
\mbox{ and} \
g \in
\mathcal{AP}(\mathbb{R},\mathbb{R}^{n}).
\end{equation}

We also introduce the following notation
 $\mathcal{PAP}
(\mathbb{R},\mathbb{R}^{n}, \mu):=\mathcal{PAP}
(\mathbb{R},\mathbb{R}^{n}, \mu,\mu)$.
\end{defn}

\begin{defn} (see \cite{mir})
 If $\mu, \nu \in \mathcal{M}$ and
  $f(t,u):\mathbb{R} \times \mathbb{R} \rightarrow \mathbb{R}^n$ is continuous, then
  $f(t,u)$ is said to be {\sl $(\mu, \nu)$-pseudo almost periodic in $t$, uniformly with respect to $u$},
  abbreviated as {\sl $(\mu, \nu)$-papu},
denoted by $f\in \mathcal{PAPU}(\mathbb{R}\times \mathbb{R},\mathbb{R}^n,\mu,\nu)$, if
  $$f=g+h, \ \mbox{ where} \ g \in \mathcal{APU}(\mathbb{R}\times \mathbb{R},\mathbb{R}^n)
  \ \mbox{ and} \
  h \in \mathcal{E}U(\mathbb{R}\times \mathbb{R},\mathbb{R}^n,\mu).$$
\end{defn}
\begin{example}
Let $\mu\in \mathcal{M}$ and
$$G(t)=[\sin(t)+\sin(\sqrt{2}t)]\cos(x)+\frac{\sin(x)}{1+t^2},\;t \in \mathbb{R}.$$
Then $G\in \mathcal{PAPU}(\mathbb{R}\times \mathbb{R},\mathbb{R},\mu)$.
\end{example}
We shall need the following two conditions: \\

 \textbf{(M.1)} For every measure $\mu \in \mathcal{M}$ and every
  $  \tau \in \mathbb{R},$
  there exist
  $\beta >0 $
  and
  a bounded interval
  $I$
  such that for every $A \in \mathcal{B},$
   $$  A \cap I = \emptyset
 \  \implies \
\mu _{\tau}(A):= \mu(\{a+\tau : \ a \in A\}) \leq \beta \mu (A).$$

\textbf{(M.2)} Measures $ \mu , \nu \in \mathcal{M} $ satisfy the following condition
$$
\lim \sup_{r \rightarrow \infty} \dfrac{\mu ([-r,r])}{\nu([-r,r])} <
\infty.
$$

\begin{lemma}(see \cite{mir}) \label{k}
Let $ \mu , \nu \in \mathcal{M} $ and suppose that conditions  (\textbf{M.1}) and (\textbf{M.2}) hold. Then
\begin{itemize}
\item  decomposition \eqref{unique} above is unique;
\item $( \mathcal{PAP} ( \mathbb{R}, \mathbb{R}^{n} , \mu, \nu), \|. \|_{\infty}) $ is a Banach space; and
\item  $\mathcal{PAP}(\mathbb{R},\mathbb{R}^{n},\mu,\nu)$ is translation invariant.
\end{itemize}
\end{lemma}
\section{Double measure pseudo almost periodic solutions}

We introduce the following notations:
$$ \sup_{t \in \mathbb{R} } \{| d_{ij} (t) | \}:= \bar{d}_{ij}, \; \; \sup_{t\in \mathbb{R}} \{| I_{i} (t) | \}: = \bar{I}_{i},$$
$$ \sup_{t \in \mathbb{R} } \{| a_{ij} (t) | \}:= \bar{a}_{ij}, \; \; \sup_{t\in \mathbb{R}} \{ | b_{ijl} (t) | \} := \bar{b}_{ijl},$$
and the following  conditions:\\

\textbf{(M.3)} $\mbox{For all} \ 1 \leq i,j,l \leq n,$
$$\{d_{ij}, a_{ij}, b_{ijl}, I_{i}\}\subset \mathcal{PAP} (\mathbb{R},\mathbb{R},\mu,\nu) .$$

\textbf{(M.4)} For all $i\in \{1,2,...n\}$,
$$[t\mapsto c_{i}(t)]\in
\mathcal{AP}(\mathbb{R},\mathbb{R}) \ \
\mbox{ and} \ \
\underset{t\in\mathbb{R}}{\inf} \{c_{i}(t)\}=c_{i}^{\ast}>0.$$

\textbf{(M.5)}
For all $p>1$ and  $1\leq j \leq n$,
$$f_j,g_j, h_j \in \mathcal{PAP} (\mathbb{R}\times\mathbb{R},\mathbb{R},\mu,\nu) $$
and
  there exist positive continuous functions
  $$L_j^{f},L_j^{g},L_j^{h} \in L^p(\mathbb{R},d\mu)\cap L^p(\mathbb{R},dx)$$
   such that
for all $t, u, v \in \mathbb{R}$,
    $$|f_j(t,u)-f_j(t,v)|<L_j^{f}(t)|u-v|,$$
    $$ |g_j(t,u)-g_j(t,v)|<L_j^{g}(t)|u-v|,$$
     $$|h_j(t,u)-h_j(t,v)|<L_j^{h}(t)|u-v|.$$

  In addition, we also assume that for $1\leq j\leq n$:
   $$f_j(t,0)=g_j(t,0)=h_j(t,0)=0,\;\mbox{for all} \ t\in\mathbb{R}.$$

\textbf{(M.6)}
$$ q_0 := \max_{ i\in \{1,2,...,n\} } \{   \frac{\sum_{j=1}^{n}\Big[ \bar{d}_{ij}\| L_{j}^{f}\|_{p} +  \bar{a}_{ij}\|L_{j}^{g}\|_{p} +  \sum_{l=1}^{n} \bar{b}_{ijl} ( \|L_{j}^{h}\|_{p} \|h_{l}\|_{\infty} + \| L_{l}^{h}\|_{p}\|h_j \|_{\infty})\Big]}{(qc^*_{i})^{\frac{1}{q}} }   \} < 1.$$

Next, define
 $$ L:= \max_{ i\in \{1,2,...,n\} } \{ \frac{\bar{I}_{i}}{c_{i}^*} \}, \;p_0: = \displaystyle \max_{ i\in \{1,2,...,n\} } \{ \frac{ \sum_{j=1}^{n} \Big[\bar{d}_{ij}\|L_{j}^{f}\|_{p}+ \bar{a}_{ij} \|L_{j}^{g }\|_{p}  +  \sum_{l=1}^{n} \bar{b}_{ijl} \|L_{j}^{h}\|_{p} \|h_{l}\|_{\infty}\Big]}{ (qc^*_{i})^{\frac{1}{q}}} \}.
$$
\begin{remark}

 If $q_0<1$, then $p_0<1$.
\end{remark}
\begin{lemma} \label{translation} Suppose that
measures
 $\mu,\nu \in \mathcal{M}$ satisfy the following requirements:
\begin{itemize}
\item $p>1$
and
condition  (\textbf{M.2}) holds;
\item $\Lambda \in \mathcal{C}(\mathbb{R}\times\mathbb{R},\mathbb{R})$ is a Lipschitz function such
that  $L^{\Lambda}\in L^p(\mathbb{R},d\mu)$; and
\item
 $y \in \mathcal{PAP} (\mathbb{R},\mathbb{R},\mu,\nu)$.
 \end{itemize}
 Then $[s\mapsto \Lambda(s,y(s-\theta))]\in \mathcal{PAP} (\mathbb{R},\mathbb{R},\mu,\nu)$, where $\theta\in \mathbb{R}$.
\end{lemma}

\begin{pf}
Since  $y \in \mathcal{PAP} (\mathbb{R},\mathbb{R},\mu,\nu)$,
it follows that
 $$y=y_{1}+y_{2},
 \ \mbox{ where} \ y_{1}\in \mathcal{AP}(\mathbb{R},\mathbb{R}) \
 \mbox{and} \
  y_{2}\in \mathcal{E}(\mathbb{R},\mathbb{R},\mu,\nu).$$
  
   Let
\begin{eqnarray*}
  \Psi(t) &=&\Lambda(t,y_{1}(t-\theta))+\Big[\Lambda\big(t,y_{1}(t-\theta)+y_{2}(t-\theta)\big)-
  \Lambda(t,y_{1}(t-\theta))\Big]
   = \Psi_{1}(t)+\Psi_{2}(t),
\end{eqnarray*}
where
$$\Psi_{1}(t)=\Lambda(t,y_{1}(t-\theta)) \ \mbox{ and} \
\Psi_{2}(t)=\Lambda\big(t,y_{1}(t-\theta)+y_{2}(t-\theta)\big)-
  \Lambda(t,y_{1}(t-\theta)).$$
  
Applying \cite{Mhamdi}, we can conclude that
$\Psi_{1}\in \mathcal{AP}(\mathbb{R},\mathbb{R})$.

Next, we  prove that
$\Psi_{2} \in \mathcal{E}(\mathbb{R},\mathbb{R},\mu,\nu)$. 
Let $z>0$, then we have:
\begin{eqnarray*}
   &&\frac{1}{\nu([-z,z])}\int_{-z}^{z}|\Psi_{2}(t)|d\mu(t)\\
    &&= \frac{1}{\nu([-z,z])}\int_{-z}^{z}
    |\Lambda\Big(t,y_{1}(t-\theta)+y_{2}(t-\theta)\Big)
    -\Lambda(t, y_{1}(t-\theta))|d\mu(t)\\
    &&\leq \frac{1}{\nu([-z,z])}\int_{-z}^{z}L^{\Lambda}(t)|y_{2}(t-\theta)|d\mu(t).
\end{eqnarray*}

    Since
    condition
    \textbf{(M.2)} holds and $y_{2}\in
\mathcal{E}(\mathbb{R},\mathbb{R},\mu,\nu)$, we get
\begin{eqnarray*}
   &&\frac{1}{\nu([-z,z])}\int_{-z}^{z}|\Psi_{2}(t)|d\mu(t)\leq\frac{1}{\nu([-z,z])}\int_{-z}^{z}L^{\Lambda}(t)|y_{2}(t-\theta)|d\mu(t)\\
    &&\leq \frac{\|y_{2}\|_{\infty}}{\nu([-z,z])}\int_{-z}^{z}L^{\Lambda}(t)d\mu(t)\\
    &&\leq \frac{\|y_{2}\|_{\infty}}{\nu([-z,z])}\Big[\int_{-z}^{z}(L^{\Lambda}(t))^pd\mu(t)\Big]^{\frac{1}{p}}
    \Big[\int_{-z}^{z}d\mu(t)\Big]^{\frac{1}{q}} \mbox{ where } \frac{1}{p}+\frac{1}{q}=1\\
   &&\leq \frac{\|y_{2}\|_{\infty}}{\nu([-z,z])^{\frac{1}{p}}}\|L^{\Lambda}\|_p \Big[\frac{\mu([-z,z])}{\nu([-z,z])}\Big]^{\frac{1}{q}}\rightarrow 0, \;as\; z\rightarrow +\infty.
    \end{eqnarray*}
    
Therefore
 $$[t\mapsto\Psi_{2}(t)]\in
\mathcal{E}(\mathbb{R},\mathbb{R},\mu,\nu) \ \mbox{ and} \ [s\mapsto \Lambda(s,y(s-\theta))]\in \mathcal{PAP} (\mathbb{R},\mathbb{R},\mu,\nu).$$

This completes the proof of Lemma~\ref{translation}.
\end{pf}\\

If  measures $\mu$ and $\nu$ are equal, then hypothesis (\textbf{M.2}) is satisfied
 and we can deduce the following corollary.
\begin{cor} \label{Ctranslation} Suppose that measure $\mu\in \mathcal{M}$ satisfies the following conditions:
\begin{itemize}
\item $p>1$;
\item $\Lambda \in \mathcal{C}(\mathbb{R}\times\mathbb{R},\mathbb{R})$ is a Lipschitz function such
that  $L^{{\Lambda}}\in L^p(\mathbb{R},d\mu)$; and
\item
 $y \in \mathcal{PAP} (\mathbb{R},\mathbb{R},\mu)$.
 \end{itemize}
 Then $[s\mapsto \Lambda(s,y(s-\theta))]\in \mathcal{PAP} (\mathbb{R},\mathbb{R},\mu)$, where $\theta\in \mathbb{R}$.
\end{cor}

\begin{lemma}\label{produit}
Let  $ \mu ,\nu\;\in \mathcal{M}$ and suppose that
 $$ y,z\in \mathcal{PAP}(\mathbb{R},
\mathbb{R},\mu,\nu).$$
Then $$y \times z\in \mathcal{PAP}(
\mathbb{R}, \mathbb{R}, \mu,\nu).$$
\end{lemma}
\begin{pf}
Since $ y,z\in \mathcal{PAP}(\mathbb{R},
\mathbb{R},\mu,\nu),$ it follows that
  $$ y = y_{1} + y_{2}  \
  \mbox{ and} \
   z= z_{1} +
z_{2}, \
\mbox{  where} \
 y_{1},z_{1}\in \mathcal{AP}(
\mathbb{R} , \mathbb{R} ) \
\mbox{ and} \
  y_{2},z_{2}\in \mathcal{E} ( \mathbb{R} , \mathbb{R} , \mu,\nu ).$$
  
    Then  $$ y
\times z= y_{1} z_{1} + y_{2} z_{1} +
y_{1} z_{2}  + y_{2} z_{2}.$$

We shall show that $y_{1} z_{1} \in \mathcal{AP}( \mathbb{R} ,
\mathbb{R} ). $ Letting   $ \varphi_{0}\in \mathcal{AP}(
\mathbb{R} , \mathbb{R} ),$ we see that:

\begin{eqnarray*}
  \| \varphi_{0}^{2} ( t  ) - \varphi_{0}^{2} (t+ \tau) \|  &=&  \|   \varphi_{0}(t) +\varphi_{0} (t + \tau )\| \cdot \| \varphi_{0} (t)-\varphi_{0} (t + \tau ) \|    \ \leq \    2  \| \varphi_{0}\|_{\infty}\cdot \varepsilon.
\end{eqnarray*}

Then $ \varphi_{0}^2\in \mathcal{AP}(
\mathbb{R} , \mathbb{R} ) $, so it follows that
$$  ( y_{1} + z_{1} )^{2} \in \mathcal{AP}(
\mathbb{R} , \mathbb{R} ) \
\mbox{  and} \
 ( y_{1} - z_{1} )^{2} \in
\mathcal{AP}( \mathbb{R} , \mathbb{R} ),$$ 

since

$$ ( y_{1} + z_{1} ) \in \mathcal{AP}( \mathbb{R} , \mathbb{R} )  \
\mbox{  and} \
 ( y_{1} - z_{1} ) \in \mathcal{AP}( \mathbb{R} , \mathbb{R} ).$$
 
  Note that
  $$y_{1} \times z_{1} = \dfrac{1}{4}  \Big(  (
y_{1} + z_{1} )^{2} - ( y_{1}  - z_{1})^{2}
\Big), $$

 so we can conclude that indeed,
 $y_{1} z_{1} \in \mathcal{AP}( \mathbb{R} ,
\mathbb{R} ) $. 

Next, we shall prove that
 $$ y_{2} z_{1} + y_{1} z_{2}  + y_{2} z_{2}  \in \mathcal{E} ( \mathbb{R} , \mathbb{R} , \mu,\nu ).$$
 
  Indeed, for $z>0$, we have:
\begin{eqnarray*}
   \dfrac{1}{\nu( [-z,z]) } \int_{-z}^{z}  |  ( y_1 z_2 + y_2 z_1 +
    y_2 z_2)(t)  | d \mu(t)
    \leq
    \dfrac{ \| y_1 \|_{\infty}}{ \nu ([ -z,z]) } \int_{-z}^{z} | z_2 (t) | d \mu(t) \\
   +  \dfrac{ \| z_{1} \|_{\infty}}{ \nu ([ -z,z]) } \int_{-z}^{z} | y_2 (t) | d
  \mu(t)
  + \dfrac{\| y_2 \|_{\infty}}{ \nu ([ -z,z])} \int_{-z}^{z} | z_{2} (t) | d \mu(t).
\end{eqnarray*}

Since $y_2,z_{2}\in \mathcal{E} ( \mathbb{R} , \mathbb{R} , \mu,\nu ) $, this completes the proof of Lemma~\ref{produit}.
\end{pf}\\

Next, we define the nonlinear operator $\Gamma$ as follows:
for  any $ \varphi =
( \varphi_{1},..., \varphi_{n} ) \in \mathcal{PAP}(\mathbb{R},
\mathbb{R}^n , \mu,\nu)$,
$$( \Gamma \circ \varphi) (t) :=x_{\varphi} (t) = ( \int_{- \infty}^{t} F_{1}(s) e^{-\int_{s}^{t}c_{1}(u)du} ds,..., \int_{-\infty}^{t}F_{n}(s)e^{-\int_{s}^{t}c_{n}(u)du} ds )^{T} ,   $$
and
$$ F_{i}(s) = \sum_{j=1}^{n} d_{ij}(s) f_{j} ( s,\varphi_{j}(s) ) + \sum_{j=1}^{n} a_{ij} (s) g_{j}(s,\varphi_{j}(s- \tau_{ij}))
$$ $$+ \sum_{j=1}^{n} \sum_{l=1}^{n} b_{i j l} (s) h_{j} (s, \varphi_{j}(s- \sigma_{ij})) h_{l} (s,s - \nu_{ij})
   + I_{i}(s).
$$
\begin{lemma}\label{gamma}
Suppose that conditions (\textbf{M.1})-(\textbf{M.6}) hold. Then $ \Gamma$  maps
$ \mathcal{PAP}(\mathbb{R}, \mathbb{R}^n , \mu ,\nu)$ into itself.
\end{lemma}
\begin{pf}
Let $ \varphi =
( \varphi_{1},..., \varphi_{n} ) \in \mathcal{PAP}(\mathbb{R},
\mathbb{R}^n , \mu,\nu)$. Then the function
 \begin{eqnarray}
F_{i} : s &\mapsto & \sum_{i=1}^{n}  d_{ij} (s) f_{j}(s,\varphi_{j} (s)) + \sum_{j=1}^{n} a_{ij}(s) g_{j}(s,\varphi_{j} (s- \tau_{ij}))  \\
& + &  \sum_{j=1}^{n}\sum_{l=1}^{n} b_{ijl}(s) h_{j}(s,s- \sigma_{ij}
) h_{l} ( \varphi_{l} (s,s - \nu_{ij} )) + I_{i}(s) ,\nonumber
\end{eqnarray}
is double measure pseudo almost periodic  for all $1 \leq i \leq n,$ by Lemmas \ref{k}, \ref{translation} and \ref{produit}.

Hence for all $1 \leq i \leq  n$, we have
$$ F_{i} = F_{i}^{1} + F_{i}^{2},    \
\mbox{where} \
 F_{i}^{1} \in
\mathcal{AP} ( \mathbb{R} , \mathbb{R} )  \
\mbox{ and} \
 F_{i}^{2} \in
\mathcal{E} ( \mathbb{R}, \mathbb{R} , \mu ,\nu). $$

 Therefore
  \begin{eqnarray}
  (\Gamma_{i} \circ\varphi) ( t) & = & \int_{- \infty }^{t} e^{-\int_{s}^{t}c_{i}(u)du} F_{i}^{1}(s) ds  + \int_{- \infty }^{t} e^{-\int_{s}^{t}c_{i}(u)du} F_{i}^{2}(s) ds \\
   &= & (\Gamma_{i}\circ F_{i}^{1})(t)  +( \Gamma_{i}\circ F_{i}^{2}) (t). \nonumber
  \end{eqnarray}  
We have to prove that $ \Gamma_{i}\circ F_{i}^{1}  \in \mathcal{AP} (
\mathbb{R} , \mathbb{R} ) , \,
 i\in \{1,2,3,...,n\}$.
To this end, note that
 \begin{eqnarray*}
 | (\Gamma_{i} \circ F_{i}^{1})(t+\tau)-( \Gamma_{i}\circ F_{i}^{1})(t)|&=&| \int_{- \infty}^{t + \tau} e^{-\int_{s}^{t+ \tau}c_{i}(u)du} F_{i}^{1}(s) d s-\int_{- \infty}^{t} e^{-\int_{s}^{t}c_{i}(u)du} F_{i}^{1}(s) d s  |\\
 &\leq & |  \int_{0}^{+ \infty} e^{-y c_{i}^{\ast}} F_{i}^{1}(t + \tau - y ) d y - \int_{0}^{+ \infty } e^{-y c_{i}^{\ast}} F_{i}^{1} ( t- y )  d y  | \\
& \leq & \int_{0}^{+ \infty} e^{- y c_{i}^{\ast}}  |  F_{i}^{1} ( t + \tau - y ) - F_{i}^{1} ( t- y ) |   d y \\
& \leq & \varepsilon  \int_{0}^{+ \infty} e^{- y c_{i}^{\ast}} dy
 \ =  \  \dfrac{\varepsilon}{c_{i}^{\ast}}.
 \end{eqnarray*}
 
 Therefore  $ \Gamma_{i} \circ F_{i}^{1} \in  \mathcal{AP} ( \mathbb{R} , \mathbb{R}^{n} ), \, 
 i\in \{1,2,3,...,n\}$.
 
 On the other hand, we can prove that $ \Gamma_{i} \circ F_{i}^{2}  \in \mathcal{E} ( \mathbb{R}, \mathbb{R} , \mu ,\nu )$, for $ 
 i\in \{1,2,3,...,n\}$.
    To this end, note that
 $$
  \int_{-z}^{z} |( \Gamma_{i}\circ F_{i}^{2}) (t)  | d\mu(t)  = \int_{-z}^{z} | \int_{-\infty}^{t} e^{-\int_{s}^{t}c_{i}(u)du} F_{i}^{2} (s) d s | d\mu(t).
 $$
 Using Fubini's Theorem, we get
 \begin{eqnarray*}
 \dfrac{1}{\nu([-z,z])} \int_{-z}^{z} | (\Gamma_{i}\circ F_{i}^{2}) (t) | d \mu(t)
 \  = \
 \dfrac{1}{\nu([-z,z])} \int_{-z}^{z} | \int_{-\infty}^{t} e^{-\int_{s}^{t}c_{i}(u)du} F_{i}^{2} (s) d s | d \mu(t) \\
  \  \leq \  \dfrac{1}{\nu([-z,z])} \int_{-z}^{z}  \int_{0}^{\infty} e^{-y c_{i}^{\ast}} | F_{i}^{2} (t - y ) |  d s d \mu(t)
  \  \leq  \
   \dfrac{1}{\nu([-z,z])} \int_{0}^{\infty}  \int_{-z}^{z} e^{-y c_{i}^{\ast}} | F_{i}^{2} (t - y ) |  d s d \mu(t),
 \end{eqnarray*}
for all $z>0$. Since
$ F_{i}^{2}  \in \mathcal{E} ( \mathbb{R},
\mathbb{R} , \mu ,\nu) ,$
 it follows by Lemma \ref{k} and
 the
dominated convergence theorem, that
$$  \Gamma_{i}\circ F_{i}^{2} \in \mathcal{E} ( \mathbb{R}, \mathbb{R} , \mu ,\nu ), \
\mbox{ for all} \
 i\in \{1,2,3,...,n\}.$$
We can thus conclude that
$$ \Gamma_{i} \circ \varphi \in
\mathcal{PAP} ( \mathbb{R} , \mathbb{R} , \mu ,\nu ) ,\
\mbox{for all}  \
i\in \{1,2,3,...,n\},$$
hence
$$
 \Gamma \circ \varphi \in \mathcal{PAP} ( \mathbb{R} ,
\mathbb{R}^{n}, \mu ,\nu  ).$$
This completes the proof of Lemma~\ref{gamma}.
\end{pf}
\begin{theorem}\label{safa}
Suppose that conditions (\textbf{M.1})--(\textbf{M.6}) hold. Then  system
\eqref{sfw}  admits  a
unique  $(\mu,\nu)-$pap solution in  $\mathbb{E}$, where
$$\mathbb{E} = \{ \psi \in \mathcal{PAP}(\mathbb{R}, \mathbb{R}^n ,
\mu,\nu): \| \psi - \varphi_{0}  \|_{\infty} \leq \frac{p_0L}{1-p_0}
\} , $$ and
$$ \varphi_{0}(t) = ( \int_{- \infty}^{t} e^{-\int_{s}^{t}c_{1}(u)du} I_{1} (s) ds,..., \int_{- \infty}^{t} e^{-\int_{s}^{t}c_{n}(u)du}I_{n}(s) d s)^T. $$
\end{theorem}
 \begin{pf}
We have
 \begin{eqnarray*}
   \| \varphi_{0}   \|_{\infty} &=&  \max_{ i\in \{1,2,...,n\} }\sup_{t \in \mathbb{R}} ( | \int_{-\infty}^{t} e^{-\int_{s}^{t}c_{i}(u)du }I_{i} (s) d s |)
    \  \leq  \  \max_{ i\in \{1,2,...,n\} }( \frac{\bar{I}_{i}}{c_{i}^{\ast}} ) := L.\end{eqnarray*}
 and
\begin{eqnarray*}
   \| \varphi \|_{\infty} &\leq & \| \varphi - \varphi_{0} \|_{\infty} + \| \varphi_{0} \|_{\infty}
     \  \leq   \  \| \varphi - \varphi_{0} \|_{\infty} +L
     \  \leq   \  \frac{L}{1-p_0}.
 \end{eqnarray*}
 Let $$ \mathbb{E} = \mathbb{E} (\varphi_{0} , p_0)= \{ \varphi \in \mathcal{PAP}(\mathbb{R}, \mathbb{R}^n , \mu ,\nu): \| \varphi - \varphi_{0} \|_{\infty} \leq \frac{p_0L}{1-p_0} \}. $$
Then for every
$ \varphi \in \mathbb{E} $, we obtain the following
\begin{eqnarray*}
  \| (\Gamma \circ \varphi) - \varphi_{0} \| _{\infty}
 &=& \max_{ i\in \{1,2,...,n\} } \sup_{t \in \mathbb{R }} \{ \Big| \displaystyle  \int_{- \infty}^{t} e^{-\int_{s}^{t}c_{i}(u)du} \sum_{j=1}^{n} [ d{ij}(s) f_{j}(s, \varphi_{j}(s)) +   a_{ij} (s) g_{j}( s,\varphi_{j}(s - \tau_{ij} )) \\
 &+&\sum_{l=1}^{n}  b_{ijl}(s) h_{j} (s, \varphi_{j} (s - \sigma_{ij} )) h_{l}( s,\varphi_{l}(s - \nu_{ij} ))] ds \Big| \} \\
 &\leq & \max_{ i\in \{1,2,...,n\} }  \sup_{t \in \mathbb{R}} \{ \int_{- \infty}^{t} e^{-\int_{s}^{t}c_{i}(u)du}  \sum_{j=1}^{n}[ \bar{d}_{ij} L_{j}^{f}(s) \| \varphi \|_{ \infty}  + \bar{a}_{ij} L_{j}^{g}(s) \| \varphi \|_{\infty} \\
 &+&  \sum_{l=1}^{n} \bar{b}_{ijl} L_{j}^{h}(s)  \|h_{l}\|_{\infty} \| \varphi \|_{\infty} ] ds \} \\
 & \leq & \max_{ i\in \{1,2,...,n\} } \{ \frac{ \sum_{j=1}^{n} [\bar{d}_{ij}\|L_{j}^{f}\|_{p}+ \bar{a}_{ij} \|L_{j}^{g }\|_{p}  +  \sum_{l=1}^{n} \bar{b}_{ijl} \|L_{j}^{h}\|_{p} \|h_{l}\|_{\infty}]}{ (qc^*_{i})^{\frac{1}{q}}} \}  \| \varphi \| _{ \infty } \\
 & \leq & p_0 \| \varphi \| _{\infty }
 \  \leq  \  p_0( \| \varphi -\varphi_0\| _{\infty } + \| \varphi_0\| _{\infty })
   \  \leq  \ \dfrac{p_0 L}{1- p_0},
 \end{eqnarray*}
 where
 $$ p_0 =  \displaystyle \max_{ i\in \{1,2,...,n\} } \{ \frac{ \sum_{j=1}^{n}\Big[ \bar{d}_{ij}\|L_{j}^{f}\|_{\infty}+ \bar{a}_{ij} \|L_{j}^{g }\|_{\infty}  + \sum_{l=1}^{n} \bar{b}_{ijl} \|L_{j}^{h}\|_{\infty} \|h_{l}\|_{\infty}\Big]}{ (qc^*_{i})^{\frac{1}{q}}} \} < 1.$$
 Therefore  $  \Gamma \circ \varphi   \in \mathbb{E}.$
Next, for all $ \phi , \psi \in \mathbb{E} ,$ we get the following
 \begin{eqnarray*}
 &&\Big|  (\Gamma_i\circ {\phi})(t) - (\Gamma_i\circ{\psi})  (t) \Big|   \leq  \displaystyle {\int_{-\infty}^{t}}   e^{-\int_{s}^{t}c_{i}(u)du} \sum_{j=1}^{n}\Big|  d_{ij} (s)\Big(f_{j}(s,\phi_j(s)) - f_j(s,\psi_{j}(s))\Big)\\
 &+&  a_{ij}(s) \Big( g_{j}(s,\phi_{j}(s- \tau_{ij})) - g_{j} (s,\psi_{j}(s- \tau_{ij} ))\Big)\\
 &+&  \sum_{l=1}^{n} b_{ijl}(s)(h_{j} (s,\phi_{j}(s- \sigma_{ij}))h_{l}(s,\phi_{l}(s- \nu_{ij}))- h_{j} (s,\psi_{j}(s- \sigma_{ij})) h_{l}(s,\psi_{l}(s-\nu_{ij} )))\Big|ds \\
 &\leq & \int_{-\infty}^{t}  e^{{-\int_{s}^{t}c_{i}(u)du}}\sum_{j=1}^{n} \Big[\bar{d}_{ij} L_{j}^{f}(s) \sup_{t \in \mathbb{R}} | \phi_{j} (t) - \psi_{j}(t) |+ \bar{a}_{ij} L_{j}^{g}(s)   \sup_{t \in \mathbb{R}}  | \phi_{j} (t) - \psi_{j}(t) |\\
 & + &   \sum_{l=1}^{n} b_{ijl}(s)| h_{j}(s,\phi_{j} (s- \sigma_{ij} ))h_{l}(s,\phi_{l}(s- \nu_{ij}))
     - h_{j}(s,\psi_{j}(s-\sigma_{ij} )) h_{l}(s,\phi_{l}(s - \nu_{ij} ))\\
      & + &  h_{j} (s,\psi_{j} (s - \sigma_{ij} (s) )) h_{l} (s,\phi_{l} (s- \nu_{ij} )) - h_{j} (s,\psi_{j} (s-\sigma_{ij} )) h_{l} (s, \psi_{l}(s- \nu_{ij} ))| \Big ] ds\\     
      & \leq & \int^{+ \infty }_{0}  e^{-c^*_{i}y}dy \sum_{j=1}^{n}\Big[ \bar{d}_{ij} L_{j}^{f}(s) \sup_{t \in \mathbb{R} } | \phi_{j}(t) - \psi_{j} (t) |
  +   \bar{a}_{ij} L_{j}^{g}(s) \sup_{t\in \mathbb{R}} | \phi_{j}(t) - \psi_{j} (t) | \\   
  & + &  \sum_{l=1}^{n} \bar{b}_{ijl}(s) ( L_{j}^{h}(s)\|h_l \|_{\infty} +  L_{l}^{h}(s)\|h_j \|_{\infty} ) \sup_{t \in \mathbb{R} } | \phi_{j}(t) - \psi_{j} (t) |\Big]ds  \\  
  & \leq &  \frac{\sum_{j=1}^{n}\Big[ \bar{d}_{ij}\| L_{j}^{f}\|_{p} +  \bar{a}_{ij}\|L_{j}^{g}\|_{p} + \sum_{l=1}^{n} \bar{b}_{ijl} ( \|L_{j}^{h}\|_{p} \|h_{l}\|_{\infty} + \| L_{l}^{h}\|_{p}\|h_j \|_{\infty})\Big]}{(qc^*_{i})^{\frac{1}{q}} } \| \phi - \psi \|_{\infty}\\
    & \leq &q_0\| \phi - \psi \|_{\infty},
 \end{eqnarray*}
where $ i= 1,..., n $. Therefore
$ \|( \Gamma \circ{\phi}) - (\Gamma \circ{\psi} )\|_{\infty}\leq q_0\| \phi - \psi \|_{\infty}.$

 Note that since $ q_0 < 1 ,$    $ \Gamma$ is a contraction and possesses a unique fixed point $ z $ which is a  $(\mu,\nu)-$pap solution of  system \eqref{sfw} in the region $ \mathbb{E}$. This completes the proof of Theorem~\ref{safa}.
  \end{pf}\\
  If the two measures $\mu$ and $\nu$ are equal, then according to the proof of  Theorem~\ref{safa}, the following corollary can be deduced.

  \begin{cor}
Suppose that conditions (\textbf{M.1}) and (\textbf{M.3})--(\textbf{M.6}) hold. Then  system
 \eqref{sfw}  admits  a
unique  $\mu-$pap solution in
$$\mathbb{E} = \{ \psi \in \mathcal{PAP}(\mathbb{R}, \mathbb{R}^n ,
\mu): \| \psi - \varphi_{0}  \|_{\infty} \leq \frac{p_0L}{1-p_0}
\} .$$

\end{cor}

  In the sequel, we shall assume that the functions $L^f_j$, $L^g_j$ and $L^h_j$ are constant. By analogy, we can prove the same results as above. In addition, by the following modifications of conditions {\bf (M.5)} and {\bf (M.6)}, the exponential stability of the solution can be obtained:\\

  \textbf{(M.7)} For all $ 1\leq j\leq n,$ there exist
 constants
 $$L_{j}^{f}, L_{j}^{g}, L_{j}^{h}, M_{j}^{f}, M_{j}^{g},M_{j}^{h}\in  \mathbb{R}^*_+$$
  such that for
all $t, x_1, x_2 \in \; \mathbb{R}$,
$$ \mid f_{j} (t,x_1) - f_{j} (t,x_2) \mid \leq L_{j}^{f} \mid x_1-x_2\mid , \quad \mid f_{j} (t,x_1) \mid \leq M_{j}^{f},$$
$$ \mid g_{j} (t,x_1) - g_{j} (t,x_2) \mid \leq L_{j}^{g} \mid x_1-x_2\mid , \quad \mid g_{j} (t,x_1) \mid \leq M_{j}^{g},$$
$$ \mid h_{j} (t,x_1) - h_{j} (t,x_2) \mid \leq L_{j}^{h} \mid x_1-x_2\mid , \quad \mid h_{j} (t,x_1) \mid \leq M_{j}^{h},$$
 $$\mbox{and } f_j(t,0)=g_j(t,0)=h_j(t,0)=0.$$

\textbf{(M.8)} There exists a nonnegative constant $q_1$
such that
$$ q_1: = \max_{ i\in \{1,2,...,n\} } \{  \frac{  \sum_{j=1}^{n}[ \bar{d}_{ij} L_{j}^{f} + \bar{a}_{ij} L_{j}^{g}
 +  \sum_{l=1}^{n} \bar{b}_{ijl} (L_{j}^{h} M_{l}^{h} + M_{j}^h L_{l}^{h} )]}{c_{i}^{\ast}}  \} < 1. $$
We let $$p_1 := \max_{ i\in \{1,2,...,n\} } \{  \frac{\sum_{j=1}^{n} [\bar{d}_{ij} L_{j}^{f} +  \bar{a}_{ij} L_{j}^{g } + \sum_{l=1}^{n} \bar{b}_{ijl} L_{j}^{h}M_{l}^{h}]}{c_{i}^{\ast}} \},
$$
and
$$ \varphi_{0}(t) := ( \int_{- \infty}^{t} e^{-\int_{s}^{t}c_{1}(u)du} I_{1} (s) ds,..., \int_{- \infty}^{t} e^{-\int_{s}^{t}c_{n}(u)du}I_{n}(s) d s)^T. $$
\begin{theorem}\label{safaa}
Suppose that conditions (\textbf{M.1})--(\textbf{M.4})
and
 (\textbf{M.7})--(\textbf{M.8}) hold. Then system \eqref{sfw}  admits a
unique $(\mu,\nu)-$pap  solution in   $\mathbb{F}$, where
$$
\mathbb{F} = \{ \psi \in \mathcal{PAP}(\mathbb{R}, \mathbb{R}^n ,
\mu,\nu): \| \psi - \varphi_{0}  \|_{\infty} \leq \frac{p_1L}{1-p_1}
\} .$$

\end{theorem}
 \begin{pf}
The following  inequality holds
\begin{eqnarray*}
  \|( \Gamma \circ{\varphi}) - \varphi_{0} \| _{\infty}  \leq  \dfrac{p_1 L}{1- p_1}.
 \end{eqnarray*}
   Therefore $ \Gamma \circ{\varphi}   \in \mathbb{F}.$
Next, for all $ \phi , \psi \in \mathbb{F} ,$
$$ \|( \Gamma\circ {\phi}) -( \Gamma \circ{\psi}) \|_{\infty}\leq q_1\| \phi - \psi \|_{\infty}.$$
  Since $ q_1< 1 ,$ it follows that
  $ \Gamma$ possesses a unique fixed point $ z $ which is a $(\mu,\nu)-$pap solution of system \eqref{sfw} in the region $ \mathbb{F}.$
  This completes the proof of Theorem~\ref{safaa}.
  \end{pf}\\
 If $\mu=\nu$, we can deduce the following result:
  \begin{cor}\label{safaaa}
Suppose that conditions (\textbf{M.1}), (\textbf{M.3})--(\textbf{M.4})
and
 (\textbf{M.7})--(\textbf{M.8}) hold. Then system \eqref{sfw}  admits a
unique $\mu-$pap  solution in
$$
\mathbb{F} = \{ \psi \in \mathcal{PAP}(\mathbb{R}, \mathbb{R}^n ,
\mu): \| \psi - \varphi_{0}  \|_{\infty} \leq \frac{p_1L}{1-p_1}
\} .$$

\end{cor}

\begin{theorem}\label{s}
     Suppose that conditions (\textbf{M.1})-(\textbf{M.4})
     and
     (\textbf{M.7})-(\textbf{M.8}) hold. Then  system \eqref{sfw} has a unique globally exponentially stable $(\mu,\nu)-$pap solution.
   \end{theorem}
   \begin{pf}
   System
    \eqref{sfw} has a unique $(\mu,\nu)-$pap solution
   $$  z(t)= ( z_{1}(t),..., z_{n} (t))^{T} \in \mathbb{E} $$ and
    $ u(t) = ( u_{1}(t),..., u_{n} (t))^{T}$ is the initial value.

    Let $ x(t)= x_{1}(t),..., x_{n}(t))^{T} $
     be an arbitrary solution of system~(\ref{sfw}) with initial value $ \varphi^{*} (t) = ( \varphi_{1}^{*}(t) ,..., \varphi_{n}^{*} (t))^{T}.$
   Let $ y_{i}(t) = x_{i}(t) -z_{i}(t) ,  \varphi_{i}(t)= \varphi_{i}^{*}(t) - u_{i}(t)$, for $ i=1...n$. Then
\begin{eqnarray} \label{sys}
y_{i}'(t)&=& -c_{i}(t) y_{i}(t) \\&+& \sum_{j=1}^{n} \Big(d_{ij} (t)( f_{j} (t,x_{j}(t))-f_{j} (t,z_{j}(t)))
  \       +  a_{ij}(t)\Big[ g_{j} (t,x_{j} (t- \tau_{ij}  ))-g_{j} (t,z_{j} (t- \tau_{ij}  ))\Big] \\
   &+&  \sum_{l=1}^{n}
  b_{ijl}(t)\Big[ h_{j}(t,x_{j}(t- \sigma_{ij}) h_{l} (t,x_{l}(t- \nu_{ij}))
\  -  \  h_{j}(t,z_{j}(t-\sigma
_{ij}))h_{l}(t,z_{l}(t-\nu _{ij}))\Big]\Big),\label{sys2}
\end{eqnarray}
where $
 i\in \{1,2,3,...,n\}$. Let $F_{i}$ be defined by
$$
F_{i}(w)= c^*_{i}-w- \sum_{j=1}^{n}\Big[\bar{d}_{ij} L_{j}^{f} +
 \bar{a}_{ij}L_{j}^{g} e^{w \tau_{ij}} +
 \sum_{l=1}^{n} \bar{b}_{ijl}( L_{j}^{h} e^{w
\sigma_{ij}} M_{l}^{h}+ M_{j}^{h}L_{l}^{h} e^{w \nu_{ij}} )\Big].
$$
 By condition  (\textbf{M.8}),
we have:
$$
F_{i}(0) = c^*_{i} - \sum_{j=1}^{n}\Big[ \bar{d}_{ij} L_{j}^{f} +
\bar{a}_{ij} L_{j}^{g} + \sum_{l=1}^{n}
\bar{b}_{ijl} (L_{j}^{h} M_{l}^{h} + M_{j}^{h} L_{l}^{h} )\Big] > 0.
$$
Thus there exists $
\varepsilon_{i}^{*} > 0 $ such that $ F_{i}(\varepsilon_{i}^{*} )=0
$ and $ F_{i}( \varepsilon_{i} > 0 )$ if $ \varepsilon_{i} \in ( 0,
\varepsilon_{i}^{*} ).$

Let $ \eta = \min \{ \varepsilon_{1}^{*} ,...,
\varepsilon_{n}^{*}  \} . $ Then
$
F_i(\eta) \geq 0$ if $i =1 ,..., n.
$
Next, there exists a nonnegative $ \lambda $ such that
$$ 0 < \lambda <
\min\{ \eta , c_{1}^{*},..., c_{n}^{*} \} \
\mbox{ and} \
 F_{i}( \lambda ) > 0 ,$$
  so for all $ i\in\{ 1,...,n\},$
\begin{equation}\label{3.7}
  \frac{1}{c^*_{i} - \lambda} \Big[\sum_{j=1}^{n} (  \bar{c}_{ij} L_{j}^{f} +  \bar{d}_{ij} L_{j}^{g} e^{ \lambda \tau_{ij}}) + \sum_{j=1}^{n} \sum_{l=1}^{n} \bar{b}_{ijl}( L_{j}^{h} e^{\lambda
\sigma_{ij}} M_{l}^{h}+ M_{j}^{h}L_{l}^{h} e^{\lambda \nu_{ij}} )\Big]< 1.
\end{equation}
Multiplying \eqref{sys}-\eqref{sys2} by $ e^{\int_{0}^{s} c_{i}(u) du} $  and integrating on
$ [0,t],$ we get
\begin{eqnarray*}
  y_{i}(t) & = & \varphi_{i} (0) e^{- \int_{0}^{t} c_{i}(u) du }+ \int_{0}^{t} e^{-\int_{s}^{t} c_{i}(u) du}  \sum_{j=1}^{n}\Big( d_{ij} (s) \Big[ f_{j}(s,y_{j}(s) + z_{j} (s)) - f_{j} (s,z_{j}(s)) \Big]  \\
   & + &  a_{ij} (s) \Big[ g_{j} (s, y_{j} (s-\tau_{ij} ) + z_{j} (s-\tau_{ij} )) - g_{j}(s, z_{j}(s- \tau_{ij} ))\Big] \\
   & + &  \sum_{l=1}^{n}
  b_{ijl}(s)\Big[ h_{j}(s,y_{j}(s- \sigma_{ij})+z_{j}(s- \sigma_{ij})) h_{l} (s,y_{j}(t- \nu_{ij})+z_{j}(t- \nu_{ij})) \\
& -& h_{j}(z_{j}(t-\sigma
_{ij}))h_{l}(z_{l}(t-\nu _{ij}))\Big] \Big)ds.
\end{eqnarray*}
Let $$M=  \displaystyle \max_{1\leq i \leq n}
\frac{c_{i}^{*}}{\sum_{j=1}^{n} \Big[(  \bar{d}_{ij} L_{j}^{f} +
\bar{ a}_{ij} L_{j}^{g} )+ \sum_{l=1}^{n}
\bar{b}_{ijl} ( L_{j}^{h} M_{l}^{h} + M_{j}^{h} L_{l}^{h} )\Big]}.$$
Clearly, $ M > 1$, and
\begin{equation*}
  \frac{1}{M}-\frac{1}{c^*_{i} - \lambda} \Big[\sum_{j=1}^{n} (  \bar{c}_{ij} L_{j}^{f} +  \bar{d}_{ij} L_{j}^{g} e^{ \lambda \tau_{ij}}) + \sum_{j=1}^{n} \sum_{l=1}^{n} \bar{b}_{ijl}( L_{j}^{h} e^{\lambda
\sigma_{ij}} M_{l}^{h}+ M_{j}^{h}L_{l}^{h} e^{\lambda \nu_{ij}} )\Big]\leq 0,
\end{equation*}
where $ 0 < \lambda < \min \{ \eta , c^*_1 , c^*_2 ,..., c^*_n
\} $ is as in \eqref{3.7}. Also,
\begin{eqnarray} \label{wit}
  \| y(t) \|_{\infty}  \leq M \| \varphi \|_{\infty} e^{-\lambda t}, \ \  t > 0.
\end{eqnarray}
 To prove inequality  (\ref{wit}), we first show that for any $ u> 1,$ the following inequality holds
\begin{equation}\label{neuf}
 \| y(t) \| _{\infty}  < u M \| \varphi \|_{\infty} e^{- \lambda t} , \  \ t > 0.
 \end{equation}
 Indeed, if \eqref{neuf} were false, there would exist some $ t_1 > 0 $ and $ i \in \{ 1,..., n\} , $ such that
 \begin{equation*}\label{dix}
   \| y(t_1) \|_{\infty} = \| y_i (t_1) \|_{\infty} = u M \| \varphi \|_{\infty} e^{- \lambda t_1}
 \end{equation*}
 and
  \begin{equation*}\label{onze}
    \| y(t) \|_{\infty} \leq u M  \| \varphi \|_{\infty} e^{-\lambda t} ,\mbox{ for every } t \in ( - \infty , t_1].
  \end{equation*}
  so we could obtain
  \begin{eqnarray*}
    | y(t_1) | &\leq & \| \varphi \|_{\infty} e^{-t_1 c^*_{i}} + \int_{0}^{t_1} e^{(- t_1 - s ) c^*_{i}}  \sum_{j=1}^{n}\Big[  \bar{d}_{ij} L_{j}^{f} \| y_{j}(s) \|_\infty +   \bar{a}_{ij} L_{j}^{g} \| y_{j} (s - \tau_{ij} ) \|_{\infty} \\
     & + &    \sum_{l=1}^{n} \bar{b}_{ijl}( L_{j}^{h}M_{l}^{h}\| y_{j} (s - \sigma_{ij} ) \|_{\infty}+ M_{j}^{h}L_{l}^{h}\| y_{j} (s - \nu_{ij} ) \|_{\infty}) \Big] ds   \\
     & \leq & \| \varphi \|_\infty e^{- t_{1} c_i^*} + \int_{0}^{t_{1}} e^{-( t_1 - s ) c^*_{i}}uM  \sum_{j=1}^{n}\Big[\bar{d}_{ij} L_{j}^{f} \| \varphi \|_\infty e^{- \lambda s} +    \bar{a}_{ij} L_{j}^{g} \| \varphi \|_{\infty} e^{-\lambda (s -\tau_{ij})}  \\
     & + &     \sum_{l=1}^{n} \bar{b}_{ijl}( L_{j}^{h}M_{l}^{h}\| \varphi \|_{\infty} e^{-\lambda (s -\sigma_{ij})}+ M_{j}^{h}L_{l}^{h}
    \| \varphi \|_{\infty} e^{-\lambda (s -\nu_{ij})} ) \Big] ds  \\
    & \leq &\| \varphi \|_\infty e^{- t_{1} c_i^*} + \int_{0}^{t_{1}} e^{-( t_1 - s ) c^*_{i}}u M\| \varphi \|_\infty e^{-\lambda s} \sum_{j=1}^{n} \Big[\bar{d}_{ij}L_{j}^{f}   +   \bar{a}_{ij} L_{j}^{g}  e^{\lambda \tau_{ij}}\\
     & + &    \sum_{l=1}^{n} \bar{b}_{ijl}( L_{j}^{h}M_{l}^{h} e^{\lambda \sigma_{ij}}+ M_{j}^{h}L_{l}^{h}
    e^{\lambda \nu_{ij}} ) \Big] ds  \\
     & \leq&  u M \| \varphi \|_\infty e^{ - \lambda t_1 }  \Big[ e^{( \lambda - a_{i*}) t_1} \Big( \frac{1}{M}-\frac{1}{c^*_{i} - \lambda} [\sum_{j=1}^{n} \{ \bar{c}_{ij} L_{j}^{f} +  \bar{d}_{ij} L_{j}^{g} e^{ \lambda \tau_{ij}}\\
      &+&  \sum_{l=1}^{n} \bar{b}_{ijl}( L_{j}^{h} e^{\lambda
\sigma_{ij}} M_{l}^{h}+ M_{j}^{h}L_{l}^{h} e^{\lambda \nu_{ij}} )\}]  \Big) \\
&+& \frac{1}{c^*_{i} - \lambda} [\sum_{j=1}^{n} (  \bar{d}_{ij} L_{j}^{f} +  \bar{a}_{ij} L_{j}^{g} e^{ \lambda \tau_{ij}})
+\sum_{j=1}^{n} \sum_{l=1}^{n} \bar{b}_{ijl}( L_{j}^{h} e^{\lambda
\sigma_{ij}} M_{l}^{h}+ M_{j}^{h}L_{l}^{h} e^{\lambda \nu_{ij}} )]\Big]\\
    & \leq &  u M  \| \varphi \|_\infty  e^{- \lambda t_1} \frac{1}{c^*_{i} - \lambda} \Big[\sum_{j=1}^{n}\Big( (  \bar{c}_{ij} L_{j}^{f} +  \bar{d}_{ij} L_{j}^{g} e^{ \lambda \tau_{ij}}) +  \sum_{l=1}^{n} \bar{b}_{ijl}( L_{j}^{h} e^{\lambda
\sigma_{ij}} M_{l}^{h}+ M_{j}^{h}L_{l}^{h} e^{\lambda \nu_{ij}} )\Big)\Big] \\
     &=& u M \| \varphi \|_\infty e^{- \lambda t_1}. \\
  \end{eqnarray*}
  Hence we could conclude that $ \| y( t_{1} )  \|_{\infty} < u M  \| \varphi \|_{\infty} e^{- \lambda t_{1}} $, which contradicts  inequality \eqref{neuf}. Note that $ u \rightarrow 1 , $
so \eqref{wit} holds. Therefore system
 \eqref{sfw}
has a unique globally exponentially stable $(\mu,\nu)-$pap solution.
This completes the proof of Theorem~\ref{s}.
   \end{pf}\\
  If $\mu=\nu$, then hypothesis (\textbf{M.2}) is satisfied and we can deduce the following corollary:

   \begin{cor}
     Suppose that conditions (\textbf{M.1}), (\textbf{M.3})-(\textbf{M.4}),
     and
     (\textbf{M.7})-(\textbf{M.8}) hold. Then  system \eqref{sfw} has a unique globally exponentially stable $\mu-$pap solution.
   \end{cor}

\section{ An application to neural networks}

Neural networks have attracted a lot of attention in recent years and especially the special case
of
the so-called
  high-order Hopfield neural networks (HOHNNs) which have
been intensively investigated by many scholars in recent years,
because of their stronger approximation characteristics, larger storage capacity, faster
convergence speed, and higher fault
tolerance than low-order Hopfeld neural networks. Many
excellent  results about their dynamic characteristics
have been obtained in e.g. \cite{Adel1,BM1,Arbi1,Hajer,Mhamdi, DMO,B1,CY}. Clearly, the study of the oscillations and dynamics of such models is an exciting new topic. 

Using the  results from of this paper, we prove the existence,  the exponential
stability, and  the uniqueness  of $(\mu,\nu)-$pap
solutions of the following models of
  high-order Hopfield
neural networks (HOHNNs) with delays:
\begin{eqnarray}\label{sfw2}
  x_{i}'(t)&=& -c_{i}(t) x_{i}(t) + \sum_{j=1}^{n} d_{ij} (t) f_{j} (t,x_{j}(t)) +
\sum_{j=1}^{n} a_{ij}(t) g_{j} (t,x_{j} (t- \tau_{ij}) ) \nonumber\\
  &+&  \sum_{j=1}^{n} \sum_{l=1}^{n} b_{ijl}(t) h_{j}(t,x_{j}(t-
\sigma_{ij})) h_{l} (t,x_{l}(t- \nu_{ij})) + I_{i}(t),
\end{eqnarray}

 where $i\in \{1,...,n\}$.

 \par
\begin{tabular}{|l|c|}
  \hline
  $n$ &  number of neurons in neural network \\
  \hline

                                                                    $ x_{i} (t)$ &  $i^{th}$ neuron at time $t$ \\
                                                                       \hline
                                                                     $f_{j}, g_{j}, h_{j}$ &  activation function of $j^{th}$ neuron \\
  \hline

                                                                   $ d_{ij}(t), a_{ij}(t), b_{ijl}(t)$  &  functions connection weights     \\

                                                                       \hline

                                                                     $I_{i}(t)$ & external inputs at time $t$  \\
                                                                       \hline

                                                                      $c_{i}(t) > 0$ & rate of $i^{th}$ neuron  \\

                                                                       \hline

                                                                      $\tau_{ij}  \geq0,  \sigma_{ij}  \geq0, \nu_{ij} \geq 0$ &
transmission delays \\
  \hline
\end{tabular}

\bigskip

 The initial conditions associated with system
  \eqref{sfw2} are of the form
$$ x_i(s)= \varphi_{i} (s) , \; s \in (- \theta , 0 ] , \; i= 1,2,...,n ,$$

In our paper we have generalized the previous results by using the notion of double measure and working with two-variable functions.\\

\textbf{Example 4.1.}
Consider the following model
\begin{eqnarray}\label{25}
  &&x_{i}'(t) = -c_{i} (t) x_{i}(t) + \sum_{j=1}^{2} d_{ij} (t) f_{j} (t,x_{j}(t))
  +
\sum_{j=1}^{2} a_{ij}(t) g_{j} (t,x_{j} (t- 1 ))\nonumber\\
   &&+ \sum_{j=1}^{2} \sum_{l=1}^{2} b_{ijl}(t) h_{j}(t,x_{j}(t-1)) h_{l} (t,x_{l}(t- 1)) + I_{i}(t), 1 \leq i \leq 2,
\end{eqnarray}
where $ c_{1}=c_{2}=2$, $g_{1}(t)= g_{2}(t)= \sin t $. Then
$$L^{g_1}=L^{g_2}=M^{g_1}=M^{g_2}=1,
\tau_{i j } = \sigma_{i j } = \nu_{i j } = 1 .$$

Measures $\mu$ and $\nu$ are  defined by the
following double weights, respectively:

\begin{equation*}
\rho _{1}(t)=e^{\sin (t)},\;t\in \mathbb{R},
\end{equation*}%
and
\begin{equation*}
\rho _{2}(t)=\left\{
\begin{array}{rl}
e^{t} & \mbox{if}\ \ t\leq 0, \\
1 & \mbox{if}\ \ t>0.%
\end{array}%
\right.
\end{equation*}%
Then we have
\begin{equation*}
\frac{2r}{e}\leq \mu([-r,r])=\int_{-r}^r e^{\sin(t)}dt\leq 2er.
\end{equation*}

We now prove that
 $\mu \in \mathcal{M}$  satisfies condition
 \textbf{(M.1)}.  
 Indeed,
  $$\sin(\tau
+a)\leq 2+\sin(a)  
\ \ \mbox{ for all} \ \
\tau \in \mathbb{R}, a\in A,$$ which implies
that $$\mu(\tau + A)\leq e^2\mu(A)
\ \ \mbox{ for all} \ \
\tau \in \mathbb{R},
$$ so by \cite{EM}, $\nu \in \mathcal{M}$ satisfies
condition
 \textbf{(M.1)}. Since
\begin{equation*}
\underset{r\rightarrow +\infty }{\mathrm{\ \limsup }}\dfrac{\mu ([-r,r])}{%
\nu ([-r,r])}=\underset{r\rightarrow +\infty }{\mathrm{\ \limsup }}\dfrac{%
\displaystyle\int_{-r}^r\rho _{1}(t)dt}{\displaystyle\int_{-r}^r\rho
_{2}(t)dt}<\infty ,
\end{equation*}
it follows that
condition  \textbf{(M.2)} 
is also satisfied. We set
$$
(d_{ij}(t))_{1 \leq i,j \leq 2} =
\begin{pmatrix} \frac{2 \sin t + e^{-t}}{10}  & \frac{\cos t}{10}\\
\frac{\sin \sqrt{2} t+  e^{-t}}{10}  &\frac{2 \cos \sqrt{2}t + e^{-t}}{10}\\
\end{pmatrix},
$$
$$  (a_{ij}(t))_{1 \leq i,j \leq 2} = \begin{pmatrix} \frac{ \cos t +  e^{-t}}{10}  & \frac{ \sin t}{10}\\
\frac{4  \cos t +  e^{-t}}{10}  & \frac{\sin t +e^{-t}}{10}
\end{pmatrix},$$

$$(I_{i}(t))_{1 \leq i,j \leq 2} = \begin{pmatrix}\frac{8 \cos \sqrt{5}t}{10}   & \\
\frac{5  \sin t + e^{-t}}{10}  &   \\
\end{pmatrix},
$$
$$
(b_{1jl})(t))_{1 \leq j,l \leq 2} = \begin{pmatrix} 0  & \frac{3 \sin \sqrt{3} t + e^{-t}}{10}  \\
0   & 0   \\
\end{pmatrix},$$

$$(b_{2jl}(t))_{1 \leq j,l \leq 2} =    \begin{pmatrix} 0   & \frac{2 \cos \sqrt{5} t+ e^{-t}}{10}  \\
0  & 0   \\
\end{pmatrix}.
$$

Therefore $$ L= \frac{4}{10},\; p_1 =\frac{75}{100} < 1 \mbox{and} \  q_1=\frac{9}{10}<1.$$ Using Theorems
\ref{safaa} and \ref{s}, we can now see that
model
\eqref{25} has a unique $(\mu,\nu)-$pap
solution which is globally exponentially stable on $$\mathbb{G}= \{ \varphi \in \mathcal{PAP}(\mathbb{R}, \mathbb{R}^n
, \mu,\nu): \| \varphi - \varphi_{0}  \|_{\infty} \leq \frac{12}{10}\}.$$

\section*{List of abbreviations}
\begin{itemize}
\item {\bf ap} - almost periodic
\item {\bf apu} - almost periodic uniformly
\item {\bf pap} - pseudo almost periodic
\item {\bf papu} - pseudo almost periodic uniformly
\item {\bf HOHNN} - high-field Hopfield neural network
\end{itemize}

\section*{Declarations}

\subsection*{Ethics approval and consent to participate}
Not applicable.
\subsection*{Consent for publication}
Not applicable.
\subsection*{Availability of data and materials}
Not applicable.
\subsection*{Competing interests}
Both authors declare that they have no competing interests.
\subsection*{Funding}
Supported by the Slovenian Research Agency grant P1-0292, N1-0114, N1-0083, N1-0064, and J1-8131.
\subsection*{Authors' contributions}
Both authors contributed equally to the writing of this paper. Both authors have read and approved the submitted manuscript.


\begin{thebibliography}{99}

\bibitem{Ait} 
E.H. Ait Dads, K. Ezzinbi, M. Miraoui,
\emph {$(\mu,\nu)-$Pseudo almost automorphic solutions for some nonautonomous differential equations,}
Int. J. Math. 26 (11) (2015) 1-21.

\bibitem{Adel1}
A.M. Alimi, C. Aouiti, F. Ch\'{e}rif, M.S. M'hamdi,
\emph {Dynamics and oscillations of generalized high-order Hopfield neural networks with mixed delays},
Neurocomputing 321 (2018) 274-295.

\bibitem{BM1} 
C. Aouiti, M.S. M'hamdi, F. Ch\'{e}rif,
\emph{The existence and the stability of weighted pseudo almost periodic solution of high-order Hopfield neural network,} 
International Conference on Artificial Neural Networks (2016) 478-485.

\bibitem{Arbi1} 
A. Arbi, F. Ch\'{e}rif, C. Aouiti, A. Touati,
\emph {Dynamics of new class of Hopfield neural networks with time-varying and distributed delays,}
Acta Math. Sci. 36 (3) (2016) 891-912.

\bibitem{EM} 
J. Blot, P. Cieutat, K. Ezzinbi,
\emph{New approach for weighted pseudo almost periodic functions under the light of measure theory, basic result and applications,}
Appl. Anal. 92 (3) (2013) 493-526.

\bibitem{Boch} 
S. Bochner,
\emph{Continuous mappings of almost automorphic and almost periodic functions,}
Proc. Natl. Acad. Sci. USA 52 (1964) 907-910.

\bibitem{Hajer} 
H. Brahmi, B. Ammar, F. Ch\'{e}rif, M.A. Alimi,
\emph{Stability and exponential synchronization of high-order Hopfield neural networks with mixed delays,}
Cybern. Syst. 48 (1) (2016), 1-21.

\bibitem{F1}
 F. Ch\'{e}rif, M. Miraoui,
\emph{ New results for a Lasota-Wazewska model},
Int. J. Biomath.  12 (2) (2019) 1950019.

\bibitem{C} 
P. Cieutat, S. Fatajou, G.M. N'Gu\'{e}r\'{e}kata,
\emph {Composition of pseudo almost periodic and pseudo almost automorphic functions and applications to evolution equations,}
Appl. Anal. 89 (1) (2010) 11-27.

\bibitem{BM}
T. Diagana,
\emph { Weighted pseudo almost periodic solution to some  differentiable equations,}
Nonlinear Anal. 68 (2008) 2250-2260.

\bibitem{mir} 
T. Diagana, K. Ezzinbi, M. Miraoui,
\emph{Pseudo almost periodic and pseudo-almost automorphic solutions to some evolution equations involving theoretical measure theory,}
CUBO 16 (2) (2014) 1-31.

\bibitem{E2015}
K. Ezzinbi, M. Miraoui,
\emph{$\mu$-Pseudo almost periodic and automorphic solutions in the $\alpha$-norm for some partial functional differential equations},
Numer. Funct. Anal. Optim. 36 (8) (2015) 991-1012.

\bibitem{E2016} 
K. Ezzinbi, M. Miraoui, A. Rebey,
\emph{ Measure pseudo-almost periodic solutions in the $\alpha$-norm to some neutral partial differential equations with delay}, 
Mediterr. J. Math. 13 (5) (2016), 3417-3431.

\bibitem{Mhamdi} 
M.S. M'hamdi, C. Aouiti, A.Touati, A.M. Alimi, V. Snasel,
\emph{Weighted pseudo almost-periodic solutions of shunting inhibitory cellular neural networks with mixed delays,}
Acta Math. Sci. 36 (6) (2016) 1662-1682.

\bibitem{mir3} 
M. Miraoui,
\emph{Existence of $\mu-$pseudo almost periodic solutions to some evolution equations,}
Math. Methods Appl. Sci. 40 (13) (2017) 4716-4726.

\bibitem{mir2} 
M. Miraoui,
\emph{Pseudo almost automorphic solutions for some differential equations with reflection of the argument,}
Numer. Funct. Anal. Optim. 38(3) (2017) 376-394.

\bibitem{mir4}  
M. Miraoui, N. Yaakoubi,
\emph{Measure pseudo almost periodic solutions of shunting inhibitory cellular neural networks with mixed delays},
Numer. Funct. Anal. Optim.  40 (5) (2019) 571-585.

\bibitem{GMN} 
G.M. N'Gu\'{e}r\'{e}kata,
\emph {Almost automorphic and almost periodic functions in abstract spaces,}
 Kluwer Academic Publishers, New-York, 2001.

 \bibitem{D2}  
 N.S. Papageorgiou, V.D. Radulescu, D.D. Repov\v{s},
\emph{Periodic solutions for a class of evolution inclusions},
 Comput. Math. Appl. 75 (2018) 3047-3065.

\bibitem{D1} 
N.S. Papageorgiou, V.D. Radulescu, D.D. Repov\v{s},
\emph{ Periodic solutions for implicit evolution inclusions},
Evol. Equ. Control Theory 8 (3) 2019.

 \bibitem{Ni}
 N.S. Papageorgiou, V.D. Radulescu, D.D. Repov\v{s},
\emph{ Nonlinear Analysis - Theory and Methods}, 
 Springer Monographs in Mathematics, Springer, Cham, 2019.

\bibitem{DMO} 
F. Qiu, B. Cui, W. Wu,
\emph {Global exponential stability of high  order recurrent neural network with time-varying delays,}
Appl. Math. Model. 33 (2009), 198-210.

 \bibitem{Ra}
V.D. Radulescu, D.D. Repov\v{s}, 
\emph{Partial Differential Equations with Variable Exponents. Variational Methods and Qualitative Analysis,}
Chapman and Hall, CRC, Taylor \& Francis Group, Boca Raton,  2015.

\bibitem{B1} 
B. Xiao, H. Meng,
\emph {Existence and exponential stability of positive almost-periodic solutions for high-order Hopfield neural networks,}
Appl. Math. Model. 33 (2009), 532-542.

\bibitem{CY}
Y. Yu, M. Cai,
 \emph{Existence and exponential stability of almost periodic solutions for high-order Hopfield neural networks,}
Math. Comput. Model. 47 (2008) 943-951.

\bibitem{Z}
 C.Y. Zhang,
\emph{Pseudo almost periodic solutions of some differential equations,}
J. Math. Anal. Appl. 151 (1994) 62-76.

\bibitem{Zhao} H.Y. Zhao,
\emph{Pseudo almost periodic solutions for a class of differential equation with delays depending on state,}
Adv. Nonlinear Anal. 9 (2020), 1251-1258.
\end{thebibliography}
\end{document}